\documentclass[final]{siamltex}
\usepackage{amsmath,amssymb}

\newtheorem{rem}[theorem]{Remark}

\newtheorem{thm}{\bf{Theorem}}[section]

\newtheorem{prop}[thm]{\bf{Proposition}}

\newcommand{\gb}{\beta}

\newcommand{\gC}{\Gamma}

\newcommand{\gd}{\delta}

\newcommand{\gev}{\varepsilon}

\newcommand{\gl}{\lambda}

\newcommand{\bbb}{\mathbb{B}}

\newcommand{\bbr}{\mathbb{R}}

\newcommand{\bbn}{\mathbb{N}}

\newcommand{\bbz}{\mathbb{Z}}

\newcommand{\bbe}{\mathbb{E}}

\newcommand{\bbf}{\mathbb{F}}

\newcommand{\bbs}{\mathbb{S}}

\newcommand{\cals}{\mathcal{S}}

\newcommand{\opc}{\operatorname{C}}

\newcommand{\opbuc}{\operatorname{BUC}}

\newcommand{\opb}{\operatorname{B}}

\newcommand{\oph}{\operatorname{H}}

\newcommand{\opid}{\operatorname{id}}

\newcommand{\ops}{\operatorname{S}}

\newcommand{\dhr}{\mathrel{\lhook\joinrel\relbar\kern-.8ex\joinrel

            \lhook\joinrel\rightarrow}}

\author{Patrick Guidotti \thanks{Department of Mathematics, University of California at Irvine,
        Irvine, CA 92697-3875 (gpatrick@math.uci.edu).}}
\title{A class of Free Boundary Problems with Onset of a new Phase}

\begin{document}
\bibliographystyle{siam}

\maketitle

\bigskip
\begin{abstract} 
A class of diffusion driven Free Boundary Problems is considered which is characterized by 
the initial onset of a phase and by an explicit kinematic condition for the evolution of the free 
boundary. By a domain fixing change of variables it naturally leads to coupled systems 
comprised of a singular parabolic initial boundary value problem and 
a Hamilton-Jacobi equation. Even though the one dimensional case has been thoroughly investigated, 
results as basic as well-posedness and regularity have so far not been obtained for its 
higher dimensional counterpart. In this paper a recently developed regularity theory for 
abstract singular parabolic Cauchy problems is utilized to obtain the first well-posedness results 
for the Free Boundary Problems under consideration. The derivation of elliptic regularity results 
for the underlying static singular problems will play an important role.
\end{abstract}



\pagestyle{myheadings}
\thispagestyle{plain}
\markboth{P. GUIDOTTI}{A CLASS OF FREE BOUNDARY PROBLEMS}

\section{Introduction} 
In this paper we consider a class of one phase Free Boundary Problems (FBP) characterized 
by the initial onset of a phase. Such Free Boundary Problems arise in the description 
of diffusion in polymers for instance. Under 
physically natural conditions, these problems lead to a formulation in which the phase 
is initially absent. This feature manifests itself mathematically in that some of the 
equations in the nonlinear system become singular (if they are written in a 
fixed reference domain). In spite of the fact that this type of problems have been 
intensively studied in the literature over an extended period of time, all but one 
publication do not deal with the singular case considered here in more than one space 
dimension. The simplifying assumption that the phase be initially non-empty is typically 
added to avoid the mathematical complications stemming from the singularity. 
The one dimensional case has, however, been thoroughly investigated (\cite{FP77a,FP77b,FP77c,
FMP86,CoEr88,Gui96c}) in the specific context of diffusion in polymers and, more generally, 
for diffusion driven FBPs (\cite{Ky59},\cite{Fri83} for instance). The methods used for the one dimensional 
case rely on the explicit use of the heat kernel to reduce the problem to the boundary. This 
approach cannot be used in higher dimensions because the singular behavior induced by the 
initial condition on the FBP can not be decoupled from the diffusion operator. This is due 
the the fact that the free boundary has a non trivial geometry in this case and to the fact 
that fundamental solutions (evolution operators) have not been studied for singular parabolic 
problems for which the singularity affects the underlying elliptic operator in an anisotropic way.
The higher dimensional problem has recently been studied in \cite{Gui99} but only in 
its simpler quasi-stationary form. \cite{Gui99} established well-posedness for the problem 
in a class of functions which needs to be carefully crafted and leads to an asymptotic 
expansion for the singularity which is valid in the corresponding topologies and has practical 
relevance (cf. \cite{CoEr88,Gui96c}). Regularity results, maximal regularity in particular, 
for singular parabolic equations play a crucial role in this paper and 
have been recently obtained by the author in \cite{G06} where a construction 
of the evolution operator is given for a wide class of singular parabolic problems. Maximal 
regularity is needed because the  full problem will be solved by reduction to a model problem and 
perturbation about it. Since the ensuing perturbation is of maximal order (both in the singularity 
and differentiation) optimal regularity results are necessary.
Previous abstract results were obtained by \cite{Gui96a,W98}, but they do not apply to the 
situation considered here because their validity does not cover the case of spaces of classical 
point-wise regularity, nor the case in which the singularity affects the equations anisotropically. 
Classical point-wise regularity is needed in the analysis proposed 
here since the singular parabolic equation is coupled to a Hamilton-Jacobi type 
equation modeling the front's dynamics. It has long been observed, as mentioned above, 
that one dimensional techniques can not be extended to the higher dimensional case. Summarizing, 
the presence of singular coefficients and the coupling to a Hamilton-Jacobi equations are two 
of the characterizing features of the problem under consideration. They make its analysis more 
difficult and delicate than that of the related but different classical Stefan-problem. This paper 
therefore offers a successful approach that fully overcomes these difficulties.\\
The unknowns of the problem are a function $u:\Omega_t\to\bbr$ defined on the open domain 
$\Omega _t\subset\bbr ^{n+1}\ni(x,y)$ and its unknown boundary $\Gamma _t$. The function $u$ 
measures the concentration of the penetrant in the case of diffusion in polymers. 
Many geometries can be chosen for the domain $\Omega _t$. Here a strip-like setting is chosen 
where the domain is bounded by a fixed lower, and an upper moving hyper-surface denoted by 
$\Gamma _0$ and $\Gamma _t$, respectively. Other configurations, like annulus-type domains 
are possible and interesting and the results obtained here would apply to those since they 
all would lead to the same local model problems. The system of equations satisfied by 
$(u,\Gamma _t)$ proposed in \cite{GP98} and generalizing \cite{AGL66,AS78} then reads
\begin{alignat}{2}\label{fbp1}
  \varepsilon u_t-\triangle _{x,y} u &=0\, ,&&\text{ in }\bigcup _{t>0}\{ t\}\times\Omega _t\, ,\\
  \label{fbp2}u&=g\, , &&\text{ on }(0,\infty)\times\Gamma _0\, ,\\
  \label{fbp3}-\partial _{\nu _t}u &=(1+\varepsilon u)V\, ,&&\text{ on }\bigcup _{t>0}\{ t\}
    \times\Gamma _t\, ,\\
  \label{fbp4}V&=(1+\delta \oph _t)u\, ,\qquad\qquad&&\text{ on }\bigcup _{t>0}\{ t\}
    \times\gC _t\, ,\\\label{fbp5}\Gamma _t\big |_{t=0}&=\Gamma _0\, ,&&\text{ at }t=0\, .
\end{alignat}
In the above equations $V$ denotes the front speed in normal (outward) direction $\nu _t$, and $\oph _t$ 
its mean curvature. The function $g>0$ describes the concentration profile of a reservoir of the 
diffusing molecule whereas the positive constants $\varepsilon, \delta$ arise in the 
nondimensionalization process. They measure the deviation from a threshold concentration of the 
penetrant and the strength of the curvature term, respectively. Condition \eqref{fbp3} models 
conservation of mass across the free boundary, whereas condition \eqref{fbp4} is a phenomenological 
law which is capable of capturing the behavior observed in experiments know as case II or anomalous 
diffusion. The curvature term naturally appears if one assumes the front speed to be proportional 
to the local average concentration in a reference ball of radius $\delta >0$ rather than to its 
point-wise value.
The system is clearly nonlinear and, even if the boundary conditions were linear, it still would be. 
Indeed, two different solutions have different domains of definitions and cannot be added. This typical 
feature of FBPs becomes apparent after the simple change of variables  
\begin{equation*}
  (\tau,\xi,\eta)=\bigl( t,x,\frac{y}{s(t,x)}\bigr)\, ,\: \hat u(\tau,\xi,\eta)=u\bigl(
  \tau,\xi,\eta s(\tau ,\xi)\bigr)
\end{equation*}
which transforms the problem to a corresponding one on a fixed domain. To do so, it is assumed that 
the fixed boundary of the unknown domain is given by 
$$\Gamma _0=\bbr ^{n-1}\times\{ 0\}$$ 
and that the unknown moving boundary can be described as the graph of a function 
$s(t,\cdot)$. Latter assumption is motivated by (initial) condition \eqref{fbp5}. In the new 
variables, the domain becomes the strip
$$
  S=\bbr ^{n-1}\times[0,1]
$$ 
with the obvious boundaries denoted by $\Gamma _j$, $j=0,1$. Using the old notation 
for the new variables the system now reads
\begin{align}\label{cfbp1}
  \varepsilon u_t-\triangle _x u -\frac{1+y^2|\nabla s|^2}{s^2}\partial _y^2u&=
  \varepsilon y\frac{\dot s}{s}\partial _yu-2y\frac{1}{s}\bigl(\nabla s\big | \partial _y\nabla u
  \bigr)&&-y\frac{s\triangle s-2|\nabla s|^2}{s^2}\partial _yu\\
  \label{cfbp2}u&=g\ &&\text{ on }\Gamma _0\\
  \label{cfbp3}-\frac{1+|\nabla s|^2}{s}\partial _yu&=(1+\varepsilon u)
  \dot s-\bigl(\nabla s\big |\nabla u\bigr)&&\text{ on }\Gamma _1\\
  \label{cfbp4}\dot s&=\sqrt{1+|\nabla s|^2}(1+\delta \oph _t)u
  &&\text{ on }\Gamma _1\\
  \label{cfbp5}s(0,\cdot)&=0
\end{align}
In \cite{Gui99} the quasi-stationary approximation ($\varepsilon =0$) in two space dimensions 
was considered in the presence or absence of the mean curvature term in two space dimensions. 
Here there is no restriction on the spatial dimension and the full evolutionary 
problem is analyzed but the curvature effects are neglected ($\delta =0$).\\
The main result of this paper is establishing the well-posedness of system 
\eqref{cfbp1}-\eqref{cfbp5}. To do so, appropriate function spaces have to be used which are able 
to capture the regularity of the solution as well as its asymptotic behavior at the origin. The choice of 
function spaces is limited by the simultaneous presence of a singular parabolic problem and a 
Hamilton-Jacobi equation. Spaces of classical regularity are better suited for the latter and thus 
reduce the freedom of choice for the parabolic problem. A compromise can be reached by using spaces 
of classical H\"older regularity in space and of singular H\"older behavior in time. 
A number of preparatory results are needed; they are formulated in the next section. In section 3 
elliptic and parabolic results for singular equations will be derived which play a crucial role in 
the existence proof given in the last section.
\section{Preliminaries and Setting}
In order to obtain existence results for system \eqref{cfbp1}-\eqref{cfbp5} a linearization procedure 
will be used in combination with maximal regularity results for the relevant linearized problems. 
The presence of the singularity complicates the analysis significantly since the relevant 
regularity results are not available and need to be derived first. The linearization procedure 
is meant to capture the leading order terms of the differential operators at the origin both 
with respect to their differentiation order and their singularity degree.\\
Of interest in this paper are classical solutions. If such a 
solution existed, one would be able to guess from \eqref{cfbp2} and \eqref{cfbp4}-\eqref{cfbp5} 
that, 
$$
 s(t,x)\approx t\, g(x)\text{ for }t\approx 0\, .
$$
This indicates that the relevant model problem which captures the leading order behavior of 
\eqref{cfbp1}-\eqref{cfbp3} has the form
\begin{alignat}{2}\label{mp1}
  \varepsilon u_t-\triangle _xu-\frac{1}{t^2c ^2(x)}\partial _y^2u&=f(t,x,y)&&\text{ in }(0,\infty)
  \times S\, ,\\\label{mp2}
  u&=g(t,x)&&\text{ on }(0,\infty)\times\Gamma _0\, ,\\\label{mp3}
  \frac{1}{t}\partial _yu&=h(t,x)&&\text{ on }(0,\infty)\times\Gamma _1\, ,
\end{alignat}
where in the situation at hand $g$ is independent of the time variable and $c=g$. 
The more general case of $c$ independent of $g$ is, however, of interest and motivates the different 
notation.
For reasons which can be guessed and will become apparent later, both parabolic 
($\varepsilon =1$) and elliptic ($\varepsilon=0$) regularity results are needed in 
the analysis of \eqref{mp1}-\eqref{mp3}. 
The main ingredients needed to derive such regularity results are vector-valued Fourier multiplier 
theorems and the use of spaces of singularly H\"older continuous functions in the time variable 
and the construction of an evolution operator for singular families of generators. 
The first are crucial in the analysis of the elliptic problem with time frozen and the 
latter allow for the quantitative characterization of the singular behavior in the origin 
($t=0$) both in the elliptic and the parabolic case (where maximal regularity results are 
needed).\\
For the sake of completeness the formulation of the relevant Fourier multiplier theorem and the 
definition of the classes of singular H\"older continuous functions needed in the analysis are 
given here. The basic observation illuminating the reason for their combined use will also be 
presented in this section.\\
Assume that $E$ is a given Banach space and that $T>0$ and $\beta\in(0,1)$ $[\cup\{1-\}]$, then the 
standard H\"older [Lipschitz] space is given by
$$
 \opc ^\beta\bigl([0,T],E\bigr):=\big\{ f\in\opc\bigl([0,T],E\bigr)\, \big |\, [u]_{\beta,[0T]} 
 :=\sup_{t\neq s}\frac{|u(t)-u(s)|}{|t-s|^\beta}<\infty\big\}
$$
with norm
$$
 \|\cdot\| _\beta=\|\cdot\| _{\infty ,[0,T]} +[\cdot]_{\beta,[0,T]}\, .
$$
If $\beta=1$ is chosen, one obtains the space of Lipschitz continuous functions. To distinguish it 
from the space of continuously differentiable functions the notational device $\beta=1-$ is used. 
This means that $\beta=1$ for all practical purposes except in the notation for the space which becomes 
$\opc ^{1-}$.
Singular counterparts are given by
\begin{equation}\label{shs}
 \opc ^\beta _\beta\bigl((0,T],E\bigr):=\big\{ f\in\opb\bigl((0,T],E\bigr)\, \big |\, 
  [t\mapsto t^\beta u(t)]\in\opc ^\beta\bigl((0,T],E\bigr)\}
\end{equation}
with weighted norm defined through
$$
  \|\cdot\| _{\beta,\beta}:=\| u\| _{\infty,(0,T]}+[(\cdot)^\beta u]_{\beta,(0,T]}\, .
$$
The symbol denoting the time interval in the notation for the norm will be dropped in the sequel 
with the understanding that the interval of definition does not contain the origin for 
singular spaces, whereas it does for regular ones. The following closed subspace 
of regular H\"older functions will also be useful
$$
  \opc ^\beta _0([0,T],E):=\big\{ f\in\opc ^\beta\bigl([0,T],E)\, \big |\, f(0)=0\big\}\, .
$$
The vector-valued Fourier multiplier theorem which is needed here can be found in 
\cite[Theorem 6.2]{Ama97}. We reproduce it here for the ease of the reader. The following 
symbol class
\begin{multline}\label{symb_class}
  \mathcal{S}^m(E_0,E_1):=\big\{ a\in\opc ^{n+1}\bigl(\bbr ^n\setminus\{ 0\},\mathcal{L}(E_0,E_1)
  \bigr)\,\big |\\ \sup _{\xi\in\bbr ^n}\big |[1+|\xi |^2]^{(m+|\alpha |)/2}\partial ^\alpha 
  a(\xi)\big |_{\mathcal{L}(E_0,E_1)}<\infty\, ,\: |\alpha|\leq n+1\big\}
\end{multline}
where $E_0$ and $E_1$ are Banach spaces and $m\in\bbz$, is instrumental in the formulation of 
the result. Given a symbol $a\in\mathcal{S} ^m(E_0,E_1)$, an operator can be associated to it by
$$
  a(D):=\mathcal{F}^{-1}a\,\mathcal{F}
$$
through conjugation with the vector-valued Fourier transform 
$$
  \mathcal{F}\in\mathcal{L}_{is}\Bigl(\mathcal{S}\bigl(\bbr ^n,E_j\bigr)\Bigr)
  \cap\mathcal{L}_{is}\Bigl(\mathcal{S}'\bigl(\bbr ^n,E_j\bigr)\Bigr)\, ,\: j=0,1\, ,
$$
where $\mathcal{S}$ and $\mathcal{S}'$ are the operator-valued Schwartz' spaces of fast decaying 
test functions and of tempered distributions, respectively. They are endowed with their natural 
topology.
\begin{thm}\label{afmt}\cite[Theorem 6.2]{Ama97}
Suppose that $\mathcal{B}\in\{ B,\overset{\circ}{B},b\}$ and that $m\in\bbr$. Then
\begin{equation*}
  \bigl[ a\mapsto a(D)\bigr]\in\mathcal{L}\Bigl(\ops ^m(E_0,E_1),\mathcal{L}\bigl(
  \mathcal{B}^s_{p,q}(\bbr ^n,E_0),\mathcal{B}^{s+m}_{p,q}(\bbr ^n,E_1)\bigr)\Bigr) 
\end{equation*}
for $s\in\bbr$ and $p,q\in[1,\infty]$.
\end{thm}

The symbol $B$ denotes the regular Besov spaces whereas $\overset{\circ\quad}{B^s_{p,q}}$ and $b^s_{p,q}$ 
denote the closures of $\cals$ and $B^{s+1}_{p,q}$ in $B^s_{p,q}$, respectively. Various equivalent 
definitions are given in \cite{Ama97}. Here only regular 
Besov spaces with $p=q=\infty$ of positive fractional order $s\in\bbr ^+\setminus\bbn$ are used in which 
case one has
$$
  B^s_{\infty,\infty}=\opbuc ^s(\bbr ^n,E)\, ,
$$
where the space on the right is the standard space of bounded and uniformly H\"older continuous 
functions given by
$$
  \opbuc ^s(\bbr ^n,E):=\big\{ f\in\opbuc ^{[s]}(\bbr ^n,E)\,\big |\, \partial ^\alpha f
  \in\opbuc ^{s-[s]}(\bbr ^n,E)\:\forall\: |\alpha |\leq [s]\,\big \}\, .
$$
A crucial observation connecting the dilation of symbols of type \eqref{symb_class} and singular H\"older 
spaces was already obtained in \cite{Gui99}. It allows one to deal with singular elliptic and, eventually, 
singular parabolic boundary value problems.
\begin{prop}\label{dom}\cite[Lemma 2.5]{Gui99}
Assume that $a\in\ops ^m(E_0,E_1)$ for some $m\in\bbz$. Then
$$
  \bigl[ t\mapsto t^m\sigma _ta\bigr]\in\opc ^{1-}_{1}\bigl((0,T],\ops ^m(E_0,E_1)\bigr)
$$
for $\sigma _t(a):=a(t\cdot)$.
\end{prop}\\
\begin{rem}\label{domreg}
A straightforward adaptation of the proof given in \cite{Gui99} of Proposition \ref{dom} also shows 
that
\begin{equation*}
  \bigl[ t\mapsto t^m\sigma _ta\bigr]\in\opc ^{1-}\bigl([0,T],\ops ^{m-1}(E_0,E_1)\bigr)
\end{equation*}
\end{rem}
The importance of this proposition lies in the fact that multiplication of singularly H\"older 
continuous functions is well-defined and continuous as a map defined in various combinations of spaces
\begin{align}\label{shm1}
  \opc ^\alpha_\alpha\times\opc ^\beta_\beta&\to\opc ^\beta_\beta\\\label{shm2}
  \opc ^\alpha _\alpha\times\opc ^\beta_0&\to\opc ^\beta _0\\\label{shm3}
  \opc ^\alpha _0\times\opc ^\beta _0&\to\opc ^\beta _{-\alpha}\\\label{shm4}
  \opc ^\alpha _0\times\opc ^\beta _\beta&\to\opc ^\beta _{\beta -\alpha}
\end{align}
for $\alpha ,\beta\in(0,1)\cup\{ 1-\}$ and $\beta\leq\alpha$. The proof is elementary and can 
be found in \cite{Gui99} along with the (natural) definition of the spaces for negative lower 
indices.\\
Finally optimal regularity results are yet another essential ingredient to deal with the full 
nonlinear problem. It provides maximal regularity results for the relevant class of singular 
parabolic problems in classes of singularly H\"older continuous functions which are, as they 
need to be, perfectly compatible with the corresponding results for singular elliptic problems 
in the same class of functions. The result needed here has been derived in \cite{G06} and is 
formulated in the next theorem. It gives conditions for the well-posedness of the singular 
abstract Cauchy problem
\begin{equation}\label{sacp}
  \dot u -A(t)u=f(t)\, ,\: t>0
\end{equation}
for a Banach space valued function $u:(0,T]\to E_0$ and a singular family of ``elliptic operators'' 
$A$ (that is, of generators of analytic semi-groups). The symbol $\mathcal{H}^-(E_0,\omega)$ denotes 
the class of generators of exponentially decaying semi-groups which are not necessarily strongly 
continuous (as introduced by \cite{S85}).
\begin{thm}\label{max_reg}\cite[Corollary 3.3]{G06}
Assume that $A$ satisfies the following assumptions 
\begin{align}\label{hypo1}
  \text{(i) }&A(t)\in\mathcal{H}^-(E_0,\omega)\, ,\: t>0\, ,\\\notag
  \text{(ii) }& \|\bigl[ A(t)-A(s)\bigr]A^{-1}(\tau)\|_{\mathcal{L}(E_0)}\leq c\frac{t-s}{t}
  \text{ and }\\\label{hypo2} &\|\bigl[ A(t)-A(s)\bigr](-A)^{-\rho}(\tau)\|_{\mathcal{L}(E_0)}
  \leq c(t-s)\, ,\\\label{hypo3}
  \text{(iii) }&\lim _{t\to 0}A^{-1}(t)=0\, ,
\end{align}
for some $\rho\in(1,2)$ and $0<\tau\leq s\leq t\leq T$. 
Let $f\in\opc ^\beta _\gamma\bigl((0,T],E_0\bigr)$ for some $\beta\in(0,1)$ and $\gamma=0,\gb$. 
Then \eqref{sacp} has a unique solution $u\in\opc ^\beta _\gamma\bigl((0,T],E_0\bigr)$ satisfying
$$
  \dot u ,Au\in\opc ^\beta _\gamma\text{ and } \| \dot u\| _{\beta,\gamma}+\| Au\| _{\beta,\gamma}
  \leq c\| f\| _{\beta,\gamma} \, 
$$
where $\|\cdot\| _{\beta,0}=\|\cdot\| _\beta$.
\end{thm}

Using a combination of elliptic and parabolic estimates it will be possible to obtain a satisfactory 
regularity theory for \eqref{mp1}-\eqref{mp3}. The latter will used in the analysis of the full 
evolution system \eqref{cfbp1}-\eqref{cfbp5} for $\delta=0$. 
\section{Elliptic and Parabolic Estimates}
In order to use the results presented in the previous section in the analysis of 
\eqref{cfbp1}-\eqref{cfbp5}, appropriate function spaces have to be chosen in which 
to work. The choice can typically be justified by balancing the regularities in the 
nonlinear equations in such a way that a fixed point argument can be applied to the 
set of equations. Whereas a variety of function spaces of Sobolev type are available 
for dealing with elliptic and parabolic problems, the fact that the system of interest 
contains a Hamilton-Jacobi equation for the evolution of $s$, which, in turn, appears 
in differentiated form as a coefficient in the equations for the evolution of $u$ 
restricts the choice to spaces of classical regularity for $s$. It is indeed impossible 
to work in Sobolev-Slobodeckii spaces for $u$ as the loss of regularity incurred in 
taking traces on the boundary cannot be made up for by \eqref{cfbp4} which does not possess 
any regularizing effect. It becomes clear that spaces of classical regularity for both 
$u$ and $s$ need to be chosen. In order to use the Fourier multiplier Theorem \ref{afmt} 
the choice is therefore reduced to spaces of bounded uniformly continuous functions in 
the $x$-variable. As far as the $y$ variable is concerned a good choice is given by the 
space of continuous functions as will become clear later. The ``base'' space $E_0$ for $u$ is 
therefore chosen as
\begin{equation}\label{bse}
  E_0:=\opbuc^{1+\alpha}\bigl(\bbr ^{n},\opc(0,1)\bigr)
\end{equation}
for $\alpha\in(0,1)$ and where $\opc(0,1)$ denotes the standard space of continuous functions on 
$[0,1]$. The space with one less regularity degree in $x$ would seem like a more natural choice 
but the choice made here makes it much easier to deal with inhomogeneous Neumann boundary 
conditions such as \eqref{cfbp3}. This point will become clear later in the analysis in Section 
\ref{pb}. 
In view of the singular non-autonomous nature of the problem, the standard procedures (\cite{Lun95}) 
to deal with inhomogeneous (nonlinear) boundary conditions cannot be utilized.\\
It follows that a proper choice of ``base'' space for the full evolutionary problem is given by
\begin{equation}\label{bsp1}
  \bbe _0:=\opc ^\beta _\beta\Bigl((0,T],\opbuc^{1+\alpha}\bigl(\bbr ^{n},\opc(0,1)\bigr)\Bigr)
\end{equation}
for $u$ and by
\begin{equation}\label{bsp2}
  \bbs :=\opc ^{1+\beta}\bigl([0,T],\opbuc ^{2+\alpha}(\bbr ^n)\bigr)\cap
  \opc ^{\beta}\bigl([0,T],\opbuc ^{3+\alpha}(\bbr ^n)\bigr)
\end{equation}
for $s$. A closer look at the underlying elliptic problem ($\varepsilon=0$) is necessary in order 
to better deal with the inhomogeneous boundary terms in \eqref{mp1}-\eqref{mp3}. The operator 
$-\triangle _x-\frac{1}{t^2g^2(x)}\partial _y^2 $ has clearly non-constant coefficients and 
multiplier Theorem \ref{afmt} is formulated only in the translation invariant case. It therefore 
needs to be shown that localization arguments apply in the operator-valued setting considered here.\\
The first step is, however, to obtain regularity results in the constant coefficient case.
\subsection{The Elliptic Case - Constant Coefficients}
Consider problem \eqref{mp1}-\eqref{mp3} and assume that $\varepsilon=0$ and that $c$ is a constant which 
can be assumed to be $1$ without loss of generality. Then taking a Fourier 
transform in the $x$-variables, one obtains a parameter dependent boundary value problem for 
$\hat u$ which can be solved explicitly. The solution has the structure
\begin{equation}\label{smp}
  \hat u(t,\xi,\cdot)=t^2\sigma _tc(\xi,\cdot)\hat f+t\sigma _tb(\xi,\cdot)\hat h+\sigma _t
  a(\xi,\cdot)\hat g
\end{equation}
where $\sigma _t$ denotes dilation by $t>0$ and the operator-valued symbols $a,b,c$ are given by
\begin{multline}\label{symbols}
  a(\xi,\cdot)=\frac{\cosh\bigl(|\xi |(1-\cdot)\bigr)}{\cosh(|\xi |)}\, ,\:
  b(\xi,\cdot)=\frac{\sinh(|\xi |\cdot)}{|\xi |\cosh(|\xi |)}\, ,\:
  c(\xi,\cdot)=\bigl( |\xi |^2 +C)^{-1}\, .
\end{multline}
The first two symbols have to be considered as multiplication operator valued with respect to 
the variable $y$, whereas the last contains the sectorial operator $C$ which is the operator 
$-\partial _y^2$ on $\opc(0,1)$ with domain 
\begin{equation}\label{C}
  \operatorname{dom}(C)=\{ u\in\opc ^2(0,1)\, |\, u(0)=0\, ,\: \partial_yu(1)=0\}=:\opc ^2_{0,0}(0,1)\, .
\end{equation}
The spaces
$$
 \opc ^k_0:=\{ u\in\opc ^k(0,1)\, |\, u(0)=0\}\, ,\: k=0,1
$$
will also be needed.
Proposition \ref{dom} can be used in combination with the multiplications \eqref{shm1}-\eqref{shm4} 
to obtain regularity results in weighted H\"older spaces provided \eqref{symbols} are in 
the appropriate symbol class.
\begin{thm}\label{symb_anal}
The symbols \eqref{symbols} satisfy
\begin{align*}
  a\in\mathcal{S}^0\bigl(\opc(0,1),&\opc(0,1)\bigr)\, ,\\ b\in\mathcal{S}^k&\bigl(
  \opc(0,1),\opc ^{1-k}_0(0,1)\bigr)\, ,\\
  &\quad c\in\mathcal{S}^j\bigl(\opc(0,1),\opc ^{2-j}_{0,0}(0,1)\bigr)  
\end{align*}
for $k=0,1$ and $j=0,1,2$. The notations for the spaces have to be interpreted so that 
boundary conditions are imposed only if enough regularity is available but otherwise 
disregarded.
\end{thm}
\begin{proof}
Take the last symbol first. It is known that the operator $C$ is 
sectorial and invertible (cf. \cite{Lun95}). It follows that
\begin{equation}\label{sec_est}
  \|(|\xi |^2+C)^{-1}\| _{\mathcal{L}(\opc,\opc ^{2-j})}\leq\, c\frac{1}{1+|\xi|^j}\, ,\: j=0,1,2  
\end{equation}
As for the derivatives, it is seen by induction that, for $\xi\neq 0$,
$$
 \partial ^\alpha(|\xi |^2+C)^{-1}=(|\xi |^2+C)^{-1}\sum _{|\beta |\leq |\alpha |}
 p_{\alpha,\beta}(\xi)(|\xi |^2+C)^{-|\gb |}
$$
for polynomials $p_{\alpha ,\beta}$ of degree at most $|\beta|$. The estimate therefore 
follows from \eqref{sec_est}. Next consider the ``boundary symbols'' $a$ and $b$. A simple 
series expansion reveals that both symbols are actually analytic functions of $\xi ^2$ (and $y$ of 
course). It follows that they are smooth for all $\xi\in\bbr ^n$ in spite of appearances. It 
is therefore legitimate to concentrate on the decay properties of the symbols without necessarily 
using the nowhere vanishing weight $(1+|\xi |^2)^{1/2}$ needed in the definition of the symbol class.
Observe that both are smooth functions of $y$ regardless of $\xi$. Consequently they define 
multiplication operators for $\opc ^{2-k}(0,1)$ and $k=0,1,2$. Also observe that, for $b$, 
the smoothing effect in $x$ is due to the decay in $\xi$. They clearly satisfy the required boundary 
conditions. They further satisfy the 
relations
\begin{gather*}
  \partial _y b(\xi,y)=a(\xi, 1-y)\, ,\: \partial _y^2 b(\xi,y)=|\xi |^2b(\xi,y)\, ,\\
  \partial _y a(\xi,y)=-|\xi |^2b(\xi,1-y)\, ,\: \partial _y^2a(\xi,y)=|\xi |^2a(\xi,y)\, ,
\end{gather*}
which make it easier to derive their mapping and symbol properties from one another. For instance, 
if it were known that $a,b\in\mathcal{S}^0(\opc,\opc)$ the first relation would imply that
$\partial _yb\in\mathcal{S}^0(\opc,\opc)$ and therefore that $b\in\mathcal{S}^0(\opc,\opc ^1)$. 
The last one, on the same assumptions, would show that $a\in\mathcal{S}^{-2}(\opc,\opc ^2)$.
More details are now given for the analysis of symbol $a$. Corresponding claims for the symbol 
$b$ can be obtained by similar calculations. The claim is that $a\in\mathcal{S}^0(C,C)$. It is 
plain that
$$
  \sup _{\xi\in\bbr ^n}\sup _{y\in[0,1]}|a(\xi,y)|<\infty\, .
$$
As for the first derivatives in $\xi$ one has
\begin{multline}\label{symb'}
  \partial _j a(\xi,y)=\frac{\xi _j}{|\xi|}\frac{1}{\cosh(|\xi|)}\bigl[ (1-y)\sinh(|\xi|(1-y))-
  \tanh(|\xi|)\cosh(|\xi |(1-y))\bigr]=\\\frac{\xi _j}{|\xi|}a(\xi,y)\underset{=:d(\xi,y)}{\underbrace{
  \bigl[ (1-y)\tanh(|\xi |(1-y))-\tanh(|\xi |)\bigr]}}\, .
\end{multline}
The maximum in $y$ is either taken on the boundary $y=0,1$ or in the interior. On the boundary it 
either vanishes or is exponentially decaying in $|\xi|$. It therefore only needs to be controlled  
in the interior to conclude its estimation. Only the term in brackets in the right-hand-side 
in the first line of \eqref{symb'} depends on $y$. Setting its first $y$-derivative equal to zero, 
one arrives at the equation
$$
  \frac{\tanh(|\xi |(1-y))}{(1-y)|\xi |}=\frac{1}{|\xi |\tanh(|\xi |)-1}
$$ 
which, for $|\xi |$ large, is asymptotic to
$$
  \frac{1}{(1-y)|\xi |}\sim \frac{1}{|\xi |}\, .
$$
Thus the maximum is located at $y\sim\frac{1}{|\xi |}$ for large $|\xi |$. As for the maximal 
value for large $|\xi |$ one has
\begin{equation}\label{aaa}
  \partial _j a(\xi,\frac{1}{|\xi |})\sim\frac{\xi _j}{|\xi|}a(\xi,\frac{1}{|\xi|})\frac{1}{|\xi |}
\end{equation}
which entails the desired estimate
$$
  \sup _{x\in\bbr ^n}\sup _{y\in[0,1]}(1+|\xi |^2)^{1/2}|\partial _ja(\xi,y)|<\infty\, .
$$
Equation \eqref{aaa} follows from
\begin{equation}\label{tanh}
  (1-\frac{1}{|\xi |})\tanh(|\xi |-1)-\tanh(|\xi |)
\end{equation}
observing that $\tanh(s)=1$ up to exponentially small terms in $s$ for $s>0$ large, that is,
$$
  |\tanh(s)-1|\leq ce^{-3s}\, s>0\, .
$$
The second derivative satisfies
$$
  \partial _i\partial _j a(\xi,y)=\frac{\xi _i\xi _j}{|\xi |^3}a(\xi,y)d(\xi,y)+
  \frac{\xi _j\xi _j}{|\xi |^2}a(\xi,y)d^2(\xi,y)+\frac{\xi _j}{|\xi |}a(\xi,y)\partial _jd(\xi,y)
$$
The first term in the right-hand-side can be seen to decay like ${1}/{|\xi |^2}$ in view of the 
explicit pre-factor and the estimate obtained for the first derivative \eqref{symb'}. The second 
and third terms add up to
$$
  \frac{\xi _i\xi _j}{|\xi |^2}a(\xi,y)\bigl[ (1-y)^2-1-2(1-y)\tanh(|\xi |)\tanh(|\xi |(1-y))+
  2\tanh^2(|\xi |)\bigr]
$$
Similar calculations as those performed for the first derivative show that the maximum is now attained at
$$
  y\sim \frac{2}{|\xi |}
$$
and it amounts to
$$
  \partial _i\partial _j a(\xi,\frac{2}{|\xi |})-\frac{\xi _i\xi _j}{|\xi |^3}a(\xi,\frac{2}{|\xi |})
  d(\xi,\frac{2}{|\xi |})\sim \frac{4}{|\xi |^2}\, .
$$
Latter follows from
$$
  (1-\frac{2}{|\xi|})^2-1-2(1-\frac{2}{|\xi|})\tanh(|\xi|-2)\tanh(|\xi|)+2\tanh ^2(|\xi|)\sim
  \frac{4}{|\xi|^2}
$$
which uses \eqref{tanh} again. In conclusion it is obtained that
$$
  \sup _{\xi\in\bbr ^n}\sup _{y\in[0,1]}(1+|\xi |^2)|\partial _i\partial _ja(\xi,y)|<\infty\, .
$$
Comparing with the proof of \cite[Lemma 2.2]{Gui99}, it is observed that the arguments are almost identical. 
The only difference are the additional factors like $\frac{\xi _j}{|\xi|}$ for the first derivative or 
$\frac{\xi _i\xi _j}{|\xi|^3}$ and $\frac{\xi _i\xi _j}{|\xi|^2}$ for the second. Terms containing 
these can always be handled through their explicit dependence on $\xi$ and the estimates from 
previous derivatives. The inductive argument 
used in \cite{Gui99} can therefore be adapted to the current situation and leads to the desired result.
\end{proof}\\
Combining these results with Theorem \ref{afmt}, Proposition \ref{dom} and \eqref{shm1}-\eqref{shm4}
the next important theorem is obtained. To simplify its formulation we introduce the abbreviated 
notation
$$
  \opc ^? _?\opbuc ^?\opc ^? _{?,?}
$$
for the function spaces
$$
  \opc ^? _?\Bigl((0,T],\opbuc ^?\bigl(\bbr ^n,\opc ^? _{?,?}(0,1)\bigr)\Bigr)
$$
where the question marks can be substituted by any of the relevant regularity and singularity 
parameters.
\begin{thm}\label{rtcc}
Let $\alpha,\beta\in(0,1)$ and $\gamma=0,\beta$. Given 
$$
  (f,g,h)\in\opc ^\beta _\gamma\opbuc^{1+\alpha}\opc\times\opc ^\beta _\gamma\opbuc ^{3+\alpha}
  \times\opc ^\beta _\gamma\opbuc ^{2+\alpha}
$$
there exists a unique solution $u$ of \eqref{mp1}-\eqref{mp3} for $c\equiv const$ which belongs 
to the space
$$
  \bbe _1:=\big\{ u\in\bbe _0\, \big |\, \partial ^\alpha _x\bigl(\frac{\partial _y}{t})^ju\in\bbe _0
  \text{ for }|\alpha |,j=0,1,2\text{ s.t. }|\alpha|+j\leq 2\big\}
$$
where
$$
  \bbe _0:=\opc ^\beta _\gamma\opbuc ^{1+\alpha}\opc\, .
$$
It is given by
$$
  u=\underset{=:A(t)^{-1}}{\underbrace{\mathcal{F}^{-1}t^2\sigma _tc\mathcal{F}}}\: f+
  \underset{=:R_N(t)}{\underbrace{\mathcal{F}^{-1}t\sigma _tb\mathcal{F}}}\: h+
  \underset{=:R_D(t)}{\underbrace{\mathcal{F}^{-1}\sigma _ta\mathcal{F}}}\: g\, .
$$
\end{thm}
\noindent Modulo the singularity, the terms in the above representation clearly point to their asymptotic 
behavior at the origin. This result needs to be extended to the non-constant coefficient case. 
This is done in the next section.
\begin{rem}
To make it easier to refer back to the above result in the case $\gamma =0$, the corresponding 
spaces will be denoted by $\bbe _{k,0}$, $k=0,1$.
\end{rem}
\subsection{The Elliptic Non-constant Coefficient Case}
In this section we adapt an abstract formulation of a localization argument proposed by Angenent in 
\cite{A99} (see also \cite{AHS94,Ama01}) for $\bbr ^n$ to cover the case of boundary value problems in 
an operator-valued context. It is used in combination with the Fourier multiplier Theorem 
\ref{afmt} to give the basic regularity results needed for the singular elliptic boundary 
value problem underlying \eqref{mp1}-\eqref{mp3} in the case of constant coefficients. 
It should be observed that localization arguments rely on perturbation results for 
differential operators by lower order terms. These are small compared to the leading order 
operator only in a qualitative sense. Perturbation results can therefore only be applied 
for the resolvents with large $\gl$ where the decay properties of the inverse to the leading 
order operator yields smallness of the lower order terms. The structure of the singular operators 
considered here can be exploited in order to avoid shifting the leading order operator and in order 
to obtain direct invertibility results. The standard argument of course shows that non-constant 
coefficient operators are sectorial whenever their constant coefficient counterpart is (see \cite{A99}). 
This remains valid here and can be used if the time variable is kept fixed.

The goal of this section is to extend the validity of Theorem \ref{rtcc} to the non-constant 
coefficient singular Boundary Value Problem
\begin{equation}\label{bvp}
  (\mathcal{A}(t),\gamma _0,\frac{1}{t}\gamma _1\partial _y):\bbe _1\to\bbe _0\times\partial\bbe _1
\end{equation}
where $\mathcal{A}(t)$ is the elliptic operator driving the evolution in \eqref{mp1} and the 
``boundary space'' $\partial\bbe _1$ is defined by
\begin{equation}\label{boundary_space}
  \partial\bbe _1=\opc ^\beta _\gamma\opbuc ^{3+\alpha}\times\opc ^\beta _\gamma
  \opbuc ^{2+\alpha}=:\partial _D\bbe _1\times\partial _N\bbe _1\, ,\: \gamma=0,\beta\, . 
\end{equation}
Observe that Theorem \ref{rtcc} simply states that the above singular Boundary Value Problem has 
a bounded inverse in the given topologies and for constant coefficients. The localization 
procedure proposed in \cite{A99} is abstract and is based on the concept of resolution. A resolution 
of a Banach space $E$ is simply a triple $(F,\varepsilon,\delta)$ where $F$ is a Banach space and 
the maps $\delta:E\to F$ and $\varepsilon:F\to E$ satisfy $\varepsilon\circ\delta=\opid _E$. The 
operator $\delta$ plays the role of the localizing operator, whereas $\gev$ resynthesizes the local 
contributions. If the space $E$ has subspaces of interest, they should correspond to ``similar 
subspaces'' of $F$ and should be left invariant by the resolution maps. The localization operator 
can be constructed as follows
$$
 \delta:\bbe _j\to \oplus _{k\in\bbz}\bbf _j^k\, ,\: u\mapsto 
 (u_k)_{k\in\bbz ^n}:=(u\varphi _{k,r})_{k\in\bbz ^n}\, ,\: j=0,1
$$
where $\bbf _j^k=\bbe _j$ for $k\in\bbz ^n$ and the sequence space $\oplus _{k\in\bbz}\bbf _j^k$ is 
endowed with the supremum norm $\sup _{k\in\bbz}\| u_k\| _{\bbf _j}$. The maps 
$\bigl(\varphi ^2_{k,r}\bigr)_{k\in\bbz ^n}$ are chosen as to form a smooth resolution of the 
identity in $\bbr ^n$ 
$$
  \sum _{k\in\bbz ^n}\varphi ^2_{k,r}\equiv 1
$$  
subordinated to a cover by cubes $Q_{k\in\bbz ^n,r}$ of fixed side length $r>0$ obtained 
by translation and dilation of the standard cube of side size $r=1$ centered at the origin. The 
support of $\varphi _{k,r}$ should contain $Q_{k,r}$ and be contained in the union of at most 
finitely many adjacent cubes. The synthesis operator is then given by
$$
 \gev:\oplus _{k\in\bbz}\bbf _j^k\to\bbe _j\, ,\: (u_k)_{k\in\bbz ^n}
 \mapsto\sum_{k\in\bbz ^n}u_k\varphi _{k,r}
$$
in which case the relation defining a resolution is clearly satisfied. Observe that the mappings 
$\gd$ and $\gev$ are continuous. The idea is to define an 
operator on $\oplus _{k\in\bbz}\bbf _1$ which parallels the Boundary Value Problem \eqref{bvp} 
and to use its inverse 
to approximate the desired inverse of \eqref{bvp} via the use of the maps $\varepsilon$ and $\delta$. 
Define the operator $(\mathcal{A}',\mathcal{B}')$ on $\oplus _{k\in\bbz}\bbf _1$ through
$$
  (\mathcal{A}^k,\mathcal{B}^k)u_k=\bigl( -\triangle _xu_k-\frac{1}{t^2c_k(x)}\partial _y^2u_k,
  \gamma _0 u_k,\frac{1}{t}\gamma _1\partial _yu_k\bigr)
$$
where $c_k|_{\supp(\varphi _{k,r})}=c|_{\supp(\varphi _{k,r})}$ and $c_k$ is otherwise 
smoothly extended without increasing its norm to the whole space (the details of the extension 
procedure can be found in \cite{AHS94,Ama01}).
The Boundary Value Problem $(\mathcal{A}^k,\mathcal{B}^k)$ can be made arbitrarily close to a 
constant coefficient one, say by substituting $c_k$ by $c_k(x_k)$ and making $r$ small. It is 
therefore invertible and enjoys the properties claimed in Theorem \ref{rtcc}. The diagonal 
operator $(\mathcal{A}',\mathcal{B}')$ is then invertible itself and maps 
$\oplus _{k\in\bbz}\bigl[\bbf _1\times\partial\bbf _1\bigr]$ to $\oplus _{k\in\bbz}\bbf _0$ 
where $\partial\bbf _1$ is defined in the obvious way. 
Denote by $\begin{bmatrix}A'(t)^{-1}& R'_D(t)& R'_N(t)\end{bmatrix}$ its inverse. Then
\begin{equation}\label{approximation_inverse}
  \varepsilon\begin{bmatrix}A'(t)^{-1}& R'_D(t) &R'_N(t)\end{bmatrix}\delta
\end{equation}
should represent an approximation to the solution operator for the singular boundary value problem 
\eqref{mp1}-\eqref{mp3}. In order to show this, it needs to be shown that \eqref{approximation_inverse}
is an approximate left and right inverse. Starting with the latter, compute
\begin{multline}\label{inv_fact}
  \begin{bmatrix}\mathcal{A}(t)\\\gamma _0\\\frac{1}{t}
  \gamma _1\partial_y\end{bmatrix}\varepsilon\begin{bmatrix}A'(t)^{-1}& R'_D(t)& R'_N(t)\end{bmatrix}
  \delta\\= \begin{bmatrix} \opid _{\bbe_0}+[\mathcal{A}(t),\varepsilon]A'(t)^{-1}\delta& 
  [\mathcal{A}(t)\, ,\varepsilon]\, R_D'(t)\delta& [\mathcal{A}(t)\, ,\varepsilon]\, R_N'(t)\delta\\ 
  0& \opid _{\gamma _0\bbe _1}& 0\\0& 0& \opid_{\gamma _1\partial _y\bbe _1}\end{bmatrix}
\end{multline}
where $[\mathcal{A}(t),\varepsilon]=\mathcal{A}(t)\varepsilon-\varepsilon \mathcal{A}(t)'$. 
Invertibility would follow if 
\begin{equation}\label{commutator}
  \opid _{\bbe _0}+[\mathcal{A}(t),\varepsilon]A'(t)^{-1}\delta
\end{equation}
were invertible. Here the specific structure of the operator $\mathcal{A}(t)$ needs to be exploited since 
standard perturbation arguments would only allow to show that the resolvent 
$\mathcal{R}(\lambda,\mathcal{A})$ exists for large enough $\gl$, in which case the lower order 
commutator term is small. The commutator term \eqref{commutator} can be computed to give
\begin{equation}\label{commutator2}
  \sum _{k\in\bbz ^n}\bigl[ 2(\nabla\varphi _{k,r}|\nabla\cdot) +\triangle\varphi _{k,r}\cdot\bigr]\bigl[
  -\triangle _x-\frac{1}{t^2c_k^2(x)}\partial _y^2\bigr]^{-1}\varphi _kf\, .
\end{equation}
To see that the inverse appearing in \eqref{commutator2} does indeed exist on the desired spaces, 
observe that $r>0$ can be chosen such that 
\begin{equation*}
  -\triangle _x-\frac{1}{t^2c_k^2(x)}\partial _y^2+\triangle _x +\frac{1}{t^2c_k^2(x_k)}\partial _y^2
  =\bigl[\frac{1}{c_k^2(x)}-\frac{1}{c_k^2(x_k)}\bigr]\frac{1}{t^2}\partial _y^2
\end{equation*}
is small in $\mathcal{L}(\bbe _1,\bbe _0)$. The realization $A^k$ of $\mathcal{A}^k$ with homogeneous 
boundary conditions is therefore invertible and the constant coefficient estimates of Theorem 
\eqref{rtcc} carry over to it (with larger constants, of course). Now, the estimates
\begin{multline*}
  \| 2\nabla\varphi _{k,r}\nabla\bigl[-\triangle _x-\frac{1}{t^2c_k^2(x)}\partial _y^2\bigr]^{-1}
  \| _{\mathcal{L}(\bbe _0)}\leq \, cT\text{ and }\\
  \| \triangle\varphi _{k,r}\bigl[-\triangle _x-\frac{1}{t^2c_k^2(x)}\partial _y^2\bigr]^{-1}
  \| _{\mathcal{L}(\bbe _0)}\leq \, cT^{2}
\end{multline*}
show that the commutator is indeed small, at least for $T$ small enough. Here the decay properties 
of the resolvent as $t\to 0$ are used as a substitute for making $\gl$ large in the standard argument. 
Heuristically this makes sense, since making $t$ small, just as making $\gl$ large, makes the operator 
``more and more elliptic''. This is due to the specific nature of the singular operator.
Since the difficulty of the problem stems from the origin, making $T$ small does not in any 
way weaken the result.\\
Next it needs to be shown that \eqref{approximation_inverse} is also a good approximation for a 
left inverse. To that end, observe that
$$
  u-\varepsilon\bigl(\mathcal{A}'(t),\mathcal{B}'(t)\bigr)^{-1}\delta\bigl(\mathcal{A}(t),
  \mathcal{B}(t)\bigr)u
$$
is the same as
\begin{equation}\label{dev_from_id}
  u-\Bigl(\varepsilon A'(t)^{-1}\delta\mathcal{A}(t)u+\varepsilon R'_D(t)\delta\gamma _0u+
  \varepsilon R'_N(t)\delta\gamma _1\frac{\partial _y}{t}u\Bigr)=\varepsilon A'(t)^{-1}
  [\mathcal{A}(t),\delta]u
\end{equation}
because
$$
 u=\varepsilon A'(t)^{-1}\mathcal{A}'(t)\delta u+\varepsilon R_D'(t)\delta\gamma _0u+
 \varepsilon R'_N(t)\gamma _1\frac{\partial _y}{t}u
$$
by definition. The notation
$$
  [\mathcal{A}(t),\delta]=\mathcal{A}'(t)\delta -\delta \mathcal{A}(t)
$$
was used in \eqref{dev_from_id}. It follows that
$$
  \varepsilon\bigl(\mathcal{A}'(t),\mathcal{B}'(t)\bigr)^{-1}\delta\bigl(\mathcal{A}(t),
  \mathcal{B}(t)\bigr)=\opid _{\bbe _1}-\varepsilon A'(t)^{-1}[\mathcal{A}(t),\delta]
$$
The term containing the commutator is lower 
order and can be estimated just like before exploiting the structure of $\mathcal{A}(t)$ in order to 
yield the invertibility of 
$$
  \opid _{\bbe _1}-\varepsilon A'(t)^{-1}[\mathcal{A}(t),\delta]
$$ 
and, thus, a left inverse for the singular boundary value problem. The result is summarized in 
the next theorem.
\begin{thm}\label{rtgc}
Let $\alpha,\beta\in(0,1)$, $\gamma=0,\beta$ and $0<c_0\leq c\in\opbuc ^{1+\alpha}$. Given 
$$
  (f,g,h)\in\opc ^\beta _\gamma\opbuc^{1+\alpha}\opc\times\opc ^\beta _\gamma\opbuc ^{3+\alpha}
  \times\opc ^\beta _\gamma\opbuc ^{2+\alpha}
$$
there exists a unique solution $u$ of \eqref{mp1}-\eqref{mp3} which belongs 
to the space
$$
  \bbe _1:=\big\{ u\in\bbe _0\, \big |\, \partial ^\alpha _x\bigl(\frac{\partial _y}{t})^ju\in\bbe _0
  \text{ for }|\alpha |,j=0,1,2\text{ s.t. }|\alpha|+j\leq 2\big\}
$$
where
$$
  \bbe _0:=\opc ^\beta _\gamma\opbuc ^{1+\alpha}\opc\, .
$$
It naturally splits into three components
$$
  u=\tilde A(t)^{-1}f+\tilde R_N(t)h+\tilde R_D(t)g
$$
with the same asymptotic behavior at the origin as in the constant coefficient case and where
\begin{gather*}
 \tilde A(t)^{-1}f=\bigl(\mathcal{A}(t),\mathcal{B}(t)\bigr) ^{-1}(f,0,0)\, ,\:
 \tilde R_D(t)g=\bigl(\mathcal{A}(t),\mathcal{B}(t)\bigr) ^{-1}(0,g,0)\\\text{ and }
 \tilde R_N(t)h=\bigl(\mathcal{A}(t),\mathcal{B}(t)\bigr) ^{-1}(0,0,h)
\end{gather*}
respectively.
\end{thm}
\begin{proof}
The only part of the proof missing is for the claim about the asymptotic behavior of the solution. It 
is obtained using the representation
\begin{equation*}
    \begin{bmatrix} T& -T[\mathcal{A}(t),\varepsilon]\, R_D'(t)\delta& -T[\mathcal{A}(t),\varepsilon]\,
    R_N'(t)\delta\\ 0& \opid _{\gamma _0\bbe _1}& 0\\0& 0& \opid_{\gamma _1\partial _y\bbe _1}
  \end{bmatrix}
\end{equation*}
for the inverse of
\begin{equation*}
  \begin{bmatrix} \opid _{\bbe _0}-[\mathcal{A}(t),\varepsilon]A'(t)^{-1}\delta& 
  [\mathcal{A},\varepsilon]R_D'(t)\delta& [\mathcal{A}(t),\varepsilon]R_N'(t)\delta\\ 
  0& \opid _{\gamma _0\bbe _1}& 0\\0& 0& 
  \opid_{\gamma _1\partial _y\bbe _1}\end{bmatrix}
\end{equation*}
where $T=\bigl[ \opid _{\bbe _0}-[\mathcal{A}(t),\varepsilon]A'(t)^{-1}\delta\bigr] ^{-1}$, 
the factorization of the resolvent implied by \eqref{inv_fact} and the mapping properties for 
the operators involved which follow from the symbol analysis combined with Theorem \ref{afmt}, 
Proposition \ref{dom} and \eqref{shm1}-\eqref{shm4}. Take for instance the entry 
$T\mathcal{A}(t)\,\varepsilon\, R_D'(t)\delta$. The claim follows from
\begin{equation*}
  T\in\mathcal{L}(\bbe _0)\, ,\: \mathcal{A}\in\mathcal{L}(\bbe _1,\bbe _0)\text{ and }
  \varepsilon R_D'\delta\in\mathcal{L}(\partial _D\bbe _1,\bbe _1)\, .
\end{equation*}
\end{proof}

From now on, the tildes on the solution operators of the non-constant coefficient case will be 
omitted. Since the corresponding operators in the constant coefficient case will no longer be used 
in the analysis, no confusion seems likely. 
\subsection{The Parabolic Problem}\label{pb}
It is now possible to return to the analysis of the model problem \eqref{mp1}-\eqref{mp3} for 
$\varepsilon=1$. If the 
boundary conditions are homogeneous, then Theorem \ref{max_reg} gives existence and regularity of a 
solution $u\in\bbe _1$ whenever $f\in\bbe _0$. The notation $\bbe _1$ (see Theorem \ref{rtcc}) has 
been used so far to describe the regularity space of the singular elliptic boundary value problem. 
In the parabolic case the same notation indicates the space
$$
 \bbe _1:=\big\{ u\in\bbe _0\, \big |\, \partial _t u\in\bbe _0\, ,\:\partial ^\alpha _x\bigl(
 \frac{\partial _y}{t})^j u\in\bbe _0\text{ for }|\alpha |,j=0,1,2\text{ s.t. }|\alpha|+j\leq 2\big\}
$$
in accordance with Theorem \ref{max_reg}.
To obtain optimal regularity results for 
$$
  A(t)=\triangle _x+\frac{1}{t^2c^2(x)}C
$$ 
it is therefore sufficient to check that conditions \eqref{hypo1}-\eqref{hypo3} are satisfied. The 
operator $C$ was defined in \eqref{C}.
Freezing coefficients arguments such as those in \cite{AHS94,Ama01} or \cite{A99} show that $A(t)$ does 
indeed generate an analytic semi-group on $E_0$ for each fixed $t>0$ which is exponentially decaying 
since $C$ is invertible (see also  the beginning paragraph of Subsection 3.2). 
Estimates \eqref{hypo2}-\eqref{hypo3} follow from the regularity theory 
developed in the previous section and the fact that $A(t)$ generates an exponentially decaying 
analytic semi-group. The latter makes it possible, in particular, to define the fractional power 
appearing in \eqref{hypo2} (see \cite{Lun95} for more details).
Condition \eqref{hypo3} follows from
$$
  \|A(t)^{-1}\| _{\mathcal{L}(\bbe _0)}\leq\, cT^2
$$
which, in its turn, follows from the observation that $(\frac{1}{c^2(x)}C)^{-1}$ is bounded in $\bbe _0$ 
and from
$$
 A(t)^{-1}=t^2\bigl( t^2\triangle _x+\frac{1}{c^2(x)}C\bigr)^{-1}\, .
$$
The first condition in \eqref{hypo2} follows from
$$
  \bigl[ A(t)-A(s)\bigr]A(\tau)^{-1}= \frac{(t^2-s^2)\tau ^2}{t^2s^2}\frac{C}{c^2(x)}
  \bigl[ -\tau ^2\triangle _{x}+\frac{C}{c^2(x)}\bigr] ^{-1}
$$
and the fact that
$$
 \|\frac{C}{c^2(x)}\bigl[ -\tau ^2\triangle _{x}+\frac{C}{c^2(x)}\bigr] ^{-1}\| 
 _{\mathcal{L}(\bbe _0)}
$$
is uniformly bounded for $\tau\in(0,T]$. The second condition in \eqref{hypo2} follows 
similarly by using abstract mapping properties of fractional powers and, in particular, that
$$
  \| \bigl[ -\tau ^2\triangle _{x}+\frac{C}{c^2(x)}\bigr] ^{-\rho +1}\| 
  _{\mathcal{L}(\bbe _0)}\leq\, c\, .
$$
It remains to be shown that problem \eqref{mp1}-\eqref{mp3} can be solved for inhomogeneous boundary 
conditions as well. The result is formulated in the next theorem.
\begin{thm}\label{iep}
Let $c\in\opbuc ^{1+\alpha}$ and
\begin{equation}\label{imax_reg}
  f\in\bbe _0\, ,\: g\in\opc ^\beta _\beta\opbuc ^{3+\alpha}\cap\opc ^{1+\beta} _\beta
  \opbuc ^{1+\alpha}\, ,h\in\opc ^\beta _\beta\opbuc ^{2+\alpha}\cap\opc ^{1+\beta} 
  _\beta\opbuc ^\alpha\, .
\end{equation}
Then there exists a unique solution $u\in\bbe _1$ of \eqref{mp1}-\eqref{mp3}.
\end{thm}
\begin{proof}
It can be assumed without loss of generality that $f\equiv 0$. Looking for a solution $u$ 
in the form
$$
  u=v+R_D(t)g+R_N(t)h
$$
one obtains that $v$ satisfies
$$
  \dot v -A(t)v=\frac{d}{dt}\bigl[ R_D(t)g+R_N(t)h\bigr]\, .
$$
Using equations \eqref{mp1}-\eqref{mp3}, it can be checked that
\begin{multline*}
  \frac{d}{dt}\bigl[ R_D(t)g\bigr]=\frac{2}{t}A(t)^{-1}\triangle _{x}R_D(t)g+R_D(t)\dot g\\
  \frac{d}{dt}\bigl[ R_N(t)h\bigr]=\frac{2}{t}A(t)^{-1}\triangle _{x}R_N(t)h+
  \frac{1}{t}R_N(t)h+R_N(t)\dot h\, .
\end{multline*}
The regularity results obtained in the previous section combined with the assumptions 
then imply that
$$
  \bigl\{ t\mapsto\frac{d}{dt}\bigl[ R_D(t)g+R_N(t)h\bigr]\bigr\}\in\bbe _0\, .
$$
Take for instance $\bigl[ t\mapsto\frac{2}{t}A(t)^{-1}\triangle _{x}R_D(t)g\bigr]$ and observe that
$$
  R_D\in\mathcal{L}(\partial _D\bbe _1,\bbe _1)\, ,\: \triangle _{x}\in\mathcal{L}(\bbe _1,\bbe _0)
  \text{ and }\bigl[t\mapsto \frac{2}{t}A(t)^{-1}\bigr]\in\mathcal{L}(\bbe _0)
$$
or $R_D(t)\dot g$ in which case the stated regularity follows from
$$
  \dot g\in\gamma _0\bbe _0\text{ and }R_D\in\mathcal{L}(\gamma _0\bbe _0,\bbe _0).
$$
The fact that $R_D\in\mathcal{L}(\gamma _0\bbe _0,\bbe _0)$ has not been explicitly proven but 
can be obtained by symbol analysis and freezing coefficients along the lines of Sections 3.1 and 3.2.
The claim then follows from Theorem \ref{max_reg}. 
\end{proof}
\section{The Hamilton-Jacobi Equation}\label{hj}
Consider the Hamilton-Jacobi equation
\begin{align}\label{hjeq1}
  s_t-\sqrt{1+|\nabla s|^2}\, v&=0 \text{ in }(0,\infty)\times\bbr ^n\, ,\\\label{hjeq2}
  s(0,\cdot)&\equiv 0\text{ on }\bbr ^n\, 
\end{align}
for a given $v\in\opc ^\beta\opbuc ^{3+\alpha}$ satisfying $v(0,\cdot)=g$. The method of 
characteristics allows one to recast this Hamilton-Jacobi equation as a system of ODEs 
in the following manner
\begin{alignat}{2}\label{hjceq1}
  \dot t &=1\, ,\: &&t(0)=0\, ,\\\label{hjceq2}
  \dot x &=-\frac{p}{\sqrt{1+|p|^2}}\, v(t,x)\, ,\: &&x(0)=\rho\in\bbr ^n\, ,\\\label{hjceq3}
  \dot z &=r-\frac{|p|^2}{\sqrt{1+|p|^2}}\, v(t,x)\, ,\: &&z(0)=0\, ,\\\label{hjceq4}
  \dot r &=\sqrt{1+|p|^2}\,\partial _t v(t,x)\, ,\: &&r(0)=g(\rho)\, \\\label{hjceq5}
  \dot p &=\sqrt{1+|p|^2}\,\nabla _xv(t,x)\, ,\: &&p(0)=0\, .
\end{alignat}
This system is easily seen to reduce to
\begin{alignat}{2}\label{rhjceq1}
  \dot x &=-\frac{p}{\sqrt{1+|p|^2}}\, v(t,x)\, ,\: &&x(0)=\rho\in\bbr ^n\, ,\\\label{rhjceq2}
  \dot p &=\sqrt{1+|p|^2}\,\nabla _xv(t,x)\, ,\: &&p(0)=0\, 
\end{alignat}
as all other unknowns can be obtained after solving this reduced system. The assumption on 
$v$ makes it possible to solve this system on a possibly small time interval $[0,T]$ which is 
independent of $\rho\in\bbr ^n$. 
Exploiting the regularity assumption on $v$ it can easily be seen that the flow mapping
$$
  (X_t,P_t):\bbr ^{2n}\to\bbr ^{2n}\, ,\: (\rho,\eta)\mapsto \bigl(x(t,\rho,\eta),p(t,\rho,\eta)\bigr)\, ,
$$
satisfies
$$
  [t\mapsto (X_t,P_t)]\in\opc ^\beta\opbuc ^{2+\alpha}(\bbr ^{2n},\bbr ^{2n})
$$
where the map is obtained by solving system \eqref{rhjceq1}-\eqref{rhjceq2} with the second initial 
condition substituted by $p(0)=\eta$. It is clearly a flow of diffeomorphisms.
Furthermore, since an equation satisfied by $(D _\rho x,D _\rho p)$ 
is easily derived from \eqref{rhjceq1}-\eqref{rhjceq2} and since 
$$
  \bigl(D _\rho x(0,\rho,\eta),D _\rho p(0,\rho,\eta)\bigr)=(\opid _{\bbr ^n},0)\, ,
$$ 
it follows that $D_\rho X_t(\rho,0)$ satisfies
$$
  \| \opid _{\bbr ^n}-D_\rho X_t(\rho,0)\|\leq\frac{1}{2}
$$
uniformly in $\rho\in\bbr ^n$. This implies that $X_t(\cdot,0)$ is a diffeomorphism and yields 
uniform estimates for the inverse of $X_t$. It follows 
that
\begin{equation}\label{diffeomreg}
  [t\mapsto (X_t)^*]\, ,\: [t\mapsto (X_t)_*]\in\opc ^\beta\mathcal{L}(\opbuc ^{2+\alpha})
\end{equation}
where we denoted the pull-back and push-forward with $X_t$ by $(X_t)^*$ and $(X_t)_*$, respectively. 
Moreover, these operators are uniformly bounded in norm in a small time interval.
\begin{thm}\label{hjeqex}
Assume that $v\in\opc ^\beta\opbuc ^{3+\alpha}$ be given with $v(0)=g\in\opbuc^{3+\alpha}$. 
Then, for $T>0$ small enough, there exists a unique solution $s$ of \eqref{hjeq1}-\eqref{hjeq2} 
on $[0,T]$ with 
$$
  s\in\opc ^{1+\beta}\opbuc ^{2+\alpha}\cap\opc ^\beta\opbuc ^{3+\alpha}\, .
$$
The existence interval is independent of $v$ in a neighborhood of $g$. Furthermore, if 
$v_1,v_2\in\opc ^\beta\opbuc ^{3+\alpha}$ are such that $v_1(0)=v_2(0)=g$, then
\begin{multline}\label{hjeqcontdep}
  \| s_1-s_2\| _{\opc ^{1+\beta}\opbuc ^{2+\alpha}\cap\opc ^\beta\opbuc ^{3+\alpha}}\\
  \leq\, cT\,\| v_1-v_2\| _{\opc ^\beta\opbuc ^{3+\alpha}}+c\,\| v_1-
  v_2\| _{\opc ^\beta\opbuc ^{2+\alpha}}\, .
\end{multline}
\end{thm}
\begin{proof}
The solution can be computed by solving the reduced system of characteristic equations 
\eqref{rhjceq1}-\eqref{rhjceq2}, which produces in particular the diffeomorphisms $X_t$. 
A minimal interval of existence can be chosen independently of the initial conditions thanks 
to the assumptions on $v$. Next
$$
  \nabla _xs(t,x)=p\bigl( t,X_t^{-1}(x),0\bigr)
$$
shows together with the discussion preceding the formulation of the theorem, in 
particular \eqref{diffeomreg}, that 
$$
  s\in\opc ^\beta\opbuc ^{3+\alpha}\, .
$$
The desired time regularity has to be backed out from equations \eqref{hjceq3}-\eqref{hjceq4}. 
The time derivative appearing in the right-hand-side of \eqref{hjceq4} is obviously not welcome 
in view of the assumptions made on $v$. It is, however, possible to rewrite the term containing 
said time derivative while integrating the equation to read
\begin{multline}\label{hjoderef}
  r(t,x)=\sqrt{1+|p(t,x)|^2}v\bigl(t,x)+\\-\int _0^t[v\nabla _x v\cdot p]
  \bigl(\tau,X_\tau(X_t^{-1}(x))\bigr)\, d\tau 
  -\int _0^t[\frac{|p|^2}{\sqrt{1+|p|^2}}v]\bigl(\tau,X_\tau(X_t^{-1}(x))\bigr)\, d\tau\, .
\end{multline}
The solution $s$ can then be obtained by one further integration from \eqref{hjceq3}. By 
\eqref{diffeomreg} and the regularity assumption on $v$ it then follows that
$$
  s\in\opc ^{1+\beta}\opbuc ^{2+\alpha}\, .
$$
The use of this argument needs to be justified. This can be done by substituting $v$ by a 
regularized version of it which is differentiable in time and for which the partial 
integration used to obtain \eqref{hjoderef} can safely be performed. Letting it converge back 
to the original function $v$ 
produces a function satisfying the modified equation. Since the solution of the original 
problem can be constructed by solving only the reduced problem \eqref{rhjceq1}-\eqref{rhjceq2} 
no trouble is encountered in taking the limit.\\
The additional estimate follows from the characteristic system \eqref{hjceq1}-\eqref{hjceq5}, 
the regularity assumptions on $v_1,v_2$ and the fact that $v_1(0)-v_2(0)=0$.
\end{proof}
\section{Existence Result}
The local existence result for regular solutions of \eqref{fbp1}-\eqref{fbp5} is based on the 
analysis of \eqref{cfbp1}-\eqref{cfbp5}. The regularity theory for linear singular elliptic 
and parabolic problems described in the previous sections and the analysis for the Hamilton-Jacobi 
equation performed in the last section will be the main tools. It should be kept in mind that 
$\varepsilon=1$ without loss of generality.
To better highlight the structure of the problem and the use of the linear regularity theory it 
is convenient to rewrite \eqref{cfbp1}-\eqref{cfbp3} using \eqref{cfbp4} in the following form
\begin{alignat}{2}\label{acfbp1}
  \dot u -\widetilde{\mathcal{A}}(t)u&=\overline{\mathcal{A}}(s)u &&\text{ in }S\\\label{acfbp2}
  u&=g&&\text{ on }\Gamma _0\\ \label{acfbp3}
  \frac{1}{t}\partial _yu&=H(s,u)&&\text{ on }\Gamma _1
\end{alignat}
where $\widetilde{\mathcal{A}}(t)=\triangle _x +\frac{1}{t^2g^2(x)}\partial _y^2+\frac{y}{t}\partial _y$ 
and the operator $\overline{\mathcal{A}}(s)$ is defined by
\begin{multline*}
  \overline{\mathcal{A}}(s)u=\bigl[\frac{1+y^2|\nabla s|^2}{s^2}-\frac{1}{t^2g^2(x)}\bigr]
  \partial _y^2u+y\bigl[\frac{\dot s}{s}-\frac{1}{t}\bigr]\partial _yu-2y\frac{1}{s}
  \bigl(\nabla s\big | \partial _y\nabla u\bigr)\\-y\frac{s\triangle s-2|\nabla s|^2}{s^2}
  \partial _yu
\end{multline*}
and $H(s,u)$ is given by
\begin{equation*}
  H(s,u)=\frac{s}{t}\frac{1}{1+|\nabla s|^2}\bigl(\nabla s\big |\nabla u\bigr)-\frac{s}{t}
  \frac{1}{\sqrt{1+|\nabla s|^2}}u(1+u)\, ,
\end{equation*}
respectively.
\begin{rem}\label{modif_gen}
It should be observed that the operator family $\widetilde{\mathcal{A}}(t)$ does not coincide with the one 
considered in the sections devoted to elliptic and parabolic regularity theory and denoted by $A(t)$. 
It is, however, easy to check that is enjoys the same ``singular'' elliptic regularity theory as a 
perturbation argument shows
$$
  \widetilde{A}(t)=A(t)\bigl[ \opid _{\bbe _1} + A^{-1}(t)\frac{y}{t}\partial _y\bigr]=
              \bigl[ \opid _{\bbe _0}+ \frac{y}{t}\partial _yA^{-1}(t)\bigr]A(t)
$$
Associated with the new family, there will be boundary solution operators $\widetilde{R}_D(t)$ 
and $\widetilde{R}_N(t)$ with corresponding mapping properties. These new boundary operators are 
actually given by
$$
  \widetilde{R}_D(t)=R_D(t)+{\widetilde A}(t)^{-1}\bigl[\frac{y}{t}\partial _yR_D(t)\bigr]\text{ and }
  \widetilde{R}_N(t)=R_N(t)+{\widetilde A}(t)^{-1}\bigl[\frac{y}{t}\partial _yR_N(t)\bigr]\, ,
$$
where $R_D$ and $R_N$ are the original boundary operators. From now on the old notation 
will apply to the new operators.
\end{rem}\\
Assume now that $s\in\bbs$, that is, that $s\in\opc ^{1+\beta}\opbuc ^{2+\alpha}\cap
\opc ^{\beta}\opbuc ^{3+\alpha}$ is chosen in a small ball about the function $[t\mapsto tg(x)]$ 
(in the natural norm of $\bbs$) and satisfying 
\begin{equation}\label{s_prop}
  s(0)=0\, , \: \dot s(0)=g\, .
\end{equation} 
It can be verified that the operator $\overline{\mathcal{A}}(s)\in\mathcal{L}(\bbe _1,\bbe _0)$ is 
small in the operator norm if the time interval length $T>0$ is chosen small (this the reason why 
the modified operator $\widetilde{A}(t)$ is introduced, see Remark \ref{modif_gen}). 
This structure will be useful in the proof of the following existence result.
\begin{thm}\label{u_ex}
Let $s\in\bbs$ be given with the above properties \eqref{s_prop}, $g\in\opbuc ^{4+\alpha}$, 
and let $T>0$ be small. 
Then there exists a unique solution $u\in\bbe _1$ of \eqref{acfbp1}-\eqref{acfbp3} which is the 
fixed-point of the operator
$$
  \Phi _1(u)=v+R_D(t)g+R_N(t)H(s,u)
$$
where $v$ solves
\begin{align}
  \dot v -A(t)v&=\overline{\mathcal{A}}(s)u-\frac{2}{t}A^{-1}(t)\bigl[\triangle _x+\frac{y}{t}
  \partial _y\bigr]\bigl[ R_D(t)g+R_N(t)H(s,u)\bigr]\label{vacfbp1}
\end{align}
Moreover, for $s_1,s_2\in\bbs$, the following estimates
\begin{multline}
  \| \gamma_1u_1-\gamma _1u_2\| _{\opc ^\beta\opbuc^{3+\beta}}\leq c\| s_1-s_2\| _{\bbs ^\beta}
  \, ,\\\| \gamma_1u_1-\gamma _1u_2\| _{\opc ^\beta\opbuc^{2+\beta}}\leq cT^\beta
  \| s_1-s_2\| _{\bbs ^\beta}
\end{multline}
hold.
\end{thm}
\begin{proof}
As already pointed, the existence proof relies on maximal regularity results obtained 
in the previous sections which are at the heart of the matter and make it now possible 
to use a simple Banach fixed-point argument. The fixed-point $u$ will be looked for in the set 
$\{ u\in\overline{\bbb} _{\bbe _1}(g,r)\, |\, u(0)=g\}$ for positive $r$ to be fixed later.
This set is endowed with the $\bbe _1$ topology and is therefore complete. It is easy but tedious 
to show that $v,R_N(t)H(s,u)$ are in $\bbe _1\cap\bbe _{0,0}$. The first regularity claim follows 
directly from Theorems \ref{rtgc} and \ref{iep} both for $v$ and $R_N(t)H(s,u)$. The vanishing property 
for $v$ follows from its regularity using the equations it satisfies. Its regularity readily implies that 
$\partial _y^2v(t=0)\equiv 0$. Since $v$ satisfies homogeneous boundary conditions 
$$
  v(y=0)=0=\partial_yv(y=1)
$$
it then is $v(t=0)\equiv 0$. As for the boundary term $R_N(t)H(s,u)$ it follows from
$$ 
  \frac{1}{t}R_N(t)\in\mathcal{L}\bigl(\opc ^\beta _0\opbuc ^{1+\alpha},\bbe _{0,0}\bigr)  
$$
combined with
$$
  tH(s,u)\in\opc ^\beta _0\opbuc ^{1+\alpha}\, ,
$$
or, alternatively, also from the equations. Finally the equations also imply that
$$
  \partial _y^2R_D(0)g\equiv 0\, ,\: R_D(0)g(y=0)=g\, ,\: \partial _yR_D(0)g(y=1)=0\, .
$$
It clearly follows that $R_D(0)g\equiv g$. Summarizing it is obtained that $\Phi _1(u)(0)=g$. 
Next a norm estimate for $\Phi _1(u)-g$ is needed. The term 
$R_D(t)g-g$ can be estimated by some fixed constant $r_1$. As for the other boundary term one has that
$$
 \| \frac{R_N(t)}{t}[tH(s,u)]\| _{\bbe _1}\leq\, cr(1+r)T^\beta\| s \| _{\bbs}
$$
and can therefore be made arbitrarily small by reducing the interval length. The above estimate is 
possible since the map $[t\mapsto tH(s,u)]$ vanishes at the origin in view of the properties of $s$. 
Next $v$ has to be estimated. Most terms can simply be estimated in a way that they produce a 
factor $T^\beta$ and can 
be made arbitrarily small by interval length reduction. This is due to the fact that
$$
  \| v\|_{\bbe _1}\leq \|\bigl( \partial _t-A(t)\bigr) ^{-1}\| _{\mathcal{L}(\bbe _0,\bbe _1)}
  \| F(s,u)\| _{\bbe _0}
$$
where $F(s,u)$ summarizes the terms in the right-hand-side of \eqref{vacfbp1}, and the specifics of 
$F$. Consider for instance the first term $\overline{\mathcal{A}}(s)u$. Since $\overline{\mathcal{A}}(s)$ 
vanishes at the origin, it can be seen that
$$
  \|\overline{\mathcal{A}}(s)\| _{\mathcal{L}(\bbe _1,\bbe _0)}\leq cT^\beta
  \| \dot s\| _{\mathcal{S}}
$$
which implies the claimed estimate. The estimate for the second term follows from the fact that
$$
  \|\frac{1}{t}A(t)^{-1}\| _{\mathcal{L}(\bbe _0)}\leq cT
$$
and the general mapping properties of the operators involved. All remaining terms but one can be 
handled similarly. All terms in $\frac{d}{dt}[tH(s,u)]$ containing either $s$ or $\nabla s$ are 
going to lead to more terms vanishing at the origin. The term
$$
  \frac{R_N(t)}{t}\bigl[\frac{\dot s}{\sqrt{1+|\nabla s|^2}}u(1+u)\bigr]
$$
behaves slightly differently and needs separate consideration since 
$\dot su(1+u)/{\sqrt{1+|\nabla s|^2}}$ does not vanish at the origin. Rewriting this term 
as
$$
  \frac{R_N(t)}{t}\bigl[\frac{\dot s}{\sqrt{1+|\nabla s|^2}}u(1+u)-g^2(1+g)\bigr]
  +\frac{R_N(t)}{t}[g^2(1+g)]
$$
and using that $u(0)=g$ up to vanishing terms of type $T^\gb$, it follows that
$$
  \| \frac{R_N(t)}{t}\bigl[\frac{\dot s}{\sqrt{1+|\nabla s|^2}}u(1+u)\bigr]\| _{\bbe _0}\leq\, 
  r_2 +cr^2\, T^\beta\| s\| _{\bbs}
$$
where the constant $r_2$ is given by
$$
  r_2=\| \bigl(\partial _t-A(t)\bigr) ^{-1}\frac{R_N(t)}{t}[g^2(1+g)]\| _{\bbe _1}\, .
$$ 
It follows that $\Phi _1$ is self-map if $r$ is chosen to satisfy $r>r_1+r_2$ and the time 
interval length is sufficiently small. It is important to observe that all estimates have 
constants which are independent of $s$ in the chosen neighborhood.\\
Next it needs to be shown that $\Phi _1$ is a contraction. The calculations are very similar to 
the above. The following terms
$$
 \Phi(u_1)-\Phi _1(u_2)=v_1-v_2 +R_N(t)\bigl[ H(s,u_1)-H(s,u_2)\bigr]
$$
need to be estimated. The difference $v_1-v_2$ satisfies
$$
  v_1-v_2=\bigl( \partial _t-A(t)\bigr) ^{-1}\bigl[ F(s,u_1)-F(s,u_2)\bigr]\, .
$$
Maximal regularity yields that
$$
  \| v_1 -v_2\| _{\bbe _1}=\|\bigl( \partial _t-A(t)\bigr) ^{-1}\bigl[ F(s,u_1)-F(s,u_2)
  \bigr]\| _{\bbe _1}\leq c\, \| F(s,u_1)-F(s,u_2)\| _{\bbe _0}
$$
so that only the last term needs to be estimated. As for the other summand, an estimate for 
$$
 \| R_N(t)\bigl[ H(s,u_1)-H(s,u_2)\bigr]\| _{\bbe _1}
$$
needs to be established. Using similar estimates as for the self-map property and 
observing that all the terms not vanishing at the origin drop out since they do 
not depend on the unknowns one arrives at
$$
  \| \Phi _1(s,u_1)-\Phi _1(s,u_2)\| _{\bbe _1}\leq\, cT^\beta\, \| u_1-u_2\| _{\bbe _1}
$$
with a constant which is independent of $s$ in the chosen neighborhood. Existence of a unique 
fixed-point follows by making the interval length short enough.\\
Next observe that any solution actually satisfies
\begin{equation}\label{decomp}
  \gamma _1u=\gamma _1R_D(t)g+\gamma _1[v+R_N(t)H(s,u)]\in\opc ^\beta\opbuc ^{3+\alpha}
\end{equation}
by Remark \ref{domreg} since $g\in\opbuc^{4+\alpha}$. This follows from 
$$
 v=\bigl(\partial _t-A(t)\bigr)^{-1}[F(s,u)+\frac{1}{t}R_N(t)g^2(1+g)]-\bigl(\partial _t-A(t)\bigr)^{-1}
   \frac{1}{t}R_N(t)g^2(1+g)
$$
and
$$
 R_N(t)H(s,u)=R_N(t)[H(s,u)-g^2(1+g)]+R_N(t)g^2(1+g)\, ,
$$
where the first terms on the right-hand-sides above belong to $\bbe _{1,0}$ and the second ones are 
not singular due to the regularity of $g$, and therefore of $g^2(1+g)$, combined with 
Remark \ref{domreg}.
The first continuous dependence result follows along the lines of the above estimates 
from the fact that the the solution to the initial boundary value problem depends smoothly 
(linearly) on the interior and boundary data $\overline{\mathcal{A}}(s)u$ and $H(s,u)$. The 
only difference stems from the fact that a stronger norm is now estimated which leads to the 
absence of the factor $T^\gb$ in front of $\| s_1-s_2\| _{\bbs}$. This factor can be 
regained if a weaker norm is estimated in view of the ``desingularizing'' properties of 
$[\partial _t-A(t)]^{-1}$, which leads to the stated continuous dependence estimate.
The estimates are based on
\begin{multline*}
  \| \gamma _1u_1 -\gamma_1 u_2\| _{\opc ^\beta\opbuc^{3+\beta}}\leq 
  \|\gamma _1\Phi _1(s_1,u_1)-\gamma _1\Phi _1(s_2,u_2)\| _{\opc ^\beta\opbuc^{3+\beta}}\\\leq
  \|\gamma _1\Phi _1(s_1,u_1)-\gamma _1\Phi _1(s_1,u_2)\| _{\opc ^\beta\opbuc^{3+\beta}}\\+
  \|\gamma _1\Phi _1(s_1,u_2)-\gamma _1\Phi _1(s_2,u_2)\| _{\opc ^\beta\opbuc^{3+\beta}}
\end{multline*}
The estimate for the first term of $\overline{\mathcal{A}}(s_1)-\overline{\mathcal{A}}(s_2)$ is typical. 
It gives 
\begin{multline*}
  \| \bigl[ \frac{t^2}{s_1^2}(1+y^2|\nabla s_1|^2)-\frac{t^2}{s_2^2}(1+y^2|\nabla s_2|^2) \bigr]
  \frac{\partial _y^2}{t^2}\| _{\mathcal{L}(\bbe _1,\bbe _{0,0})}\\\leq
  \| \frac{t^2}{s_1^2}(1+y^2|\nabla s_1|^2)-\frac{t^2}{s_2^2}(1+y^2|\nabla s_2|^2)
  \| _{\opc ^\beta\opbuc ^{1+\alpha}\opc}\, \|\frac{\partial _y^2}{t^2}\| 
  _{\mathcal{L}(\bbe _1,\bbe _0)}   
\end{multline*}
which easily entails the claim since
\begin{multline*}
  \| \frac{t^2}{s_1^2}(1+y^2|\nabla s_1|^2)-\frac{t^2}{s_2^2}(1+y^2|\nabla s_2|^2
  \| _{\opc ^\beta\opbuc ^{1+\alpha}\opc}\\\leq c\bigl(\|\dot s_1-\dot s_2\| 
  _{\opc ^\beta\opbuc ^{1+\alpha}}+\| s_1-s_2\| _{\opc ^\beta\opbuc ^{2+\alpha}}\bigr)
  \leq c\| s_1-s_2\| _\mathcal{S}\, .
\end{multline*}
One is eventually lead to
$$
 \|\gamma _1 u_1-\gamma _1u_2\| _{\opc ^\beta\opbuc^{3+\alpha}}\leq \, 
 cT^\beta\| u_1-u_2\| _{\bbe _1}+\| s_1-s_2\| _{\mathcal{S}^\gb}
$$
from which the claim follows since the first term on the right can be absorbed in the left-hand-side.
\end{proof}\\
All pieces are now in place in order to show existence for the original system of equations 
\eqref{cfbp1}-\eqref{cfbp5}. By denoting with $u=\Phi _1(s)$ the solution of \eqref{cfbp1}-\eqref{cfbp3} 
for a given $s\in\mathcal{S}$, it follows that
$$
  u\in\bbe _1\text{ and }\gamma _1u\in\opc ^\beta\opbuc ^{3+\alpha}\, .
$$
Decomposition \eqref{decomp} shows that $u$ is in a given neighborhood of $g$ for 
any choice of $s\in\mathcal{S}$ in a neighborhood of $[t\mapsto tg(x)]$. Similarly let 
$\Phi _2(u)$ be the solution of \eqref{cfbp4}-\eqref{cfbp5} constructed in section \ref{hj}. 
In this case, if $u\in\bbe _1$ with $\gamma _1u\in\opc ^\beta\opbuc ^{3+\alpha}$, 
$\Phi_2(u)\in\mathcal{S}$ will be in a neighborhood of $tg(x)$ for $u$ in a neighborhood of $g$.
It follows that a solution to the full problem can be found by producing a fixed-point $s$ for 
the map $\Phi:=\Phi _2\circ\Phi _1$ and defining $u:=\Phi _1(s)$. Summarizing
\begin{thm}\label{mth}
For any given $g\in\opbuc ^{4+\alpha}$ such that $g(x)\geq g_0>0\, ,\: x\in\bbr ^n$, there 
exists a unique local solution $(u,s)$ of the free boundary problem \eqref{cfbp1}-\eqref{cfbp5} 
such that
\begin{multline*}
  u,\dot u, \partial _x^\alpha\frac{1}{t^k}\partial _y^ku\in\opc ^\beta _\beta\opbuc ^{1+\alpha}\opc
  \text{ for }0\leq |\alpha|+k\leq 2\\
  \gamma _1u\in\opc ^\beta\opbuc ^{3+\alpha}\text{ and }\\
  s\in\opc ^{1+\beta}\opbuc ^{2+\alpha}\cap\opc ^\beta\opbuc ^{3+\alpha}\, .
\end{multline*}
Furthermore 
$$
  u=R_D(t)g+R_N(t)H(s,u)+v
$$ 
where $v$ is a solution of \eqref{vacfbp1} and the different terms have different 
asymptotic behavior in the origin. They behave 
respectively like $t^k$ for $k=0,1,2$ modulo the singular behavior built-in in the space chosen.
\end{thm}
\begin{proof}
Making the time interval as small as needed and choosing $r>0$ large enough it can be seen that
$(\Phi _1,\Phi _2)$ is a self-map on the complete set
$$
  \{u\in\overline{\bbb} _{\bbe _1}(0,r)\, |\, u(0)=g\}\times
  \{s\in\overline{\bbb} _{\bbs}(gt,g_0/2)\, |\, s(0)=0\, ,\: \dot s (0)=g\}\, .
$$
where the additional requirements $u(0)\equiv g\, ,\: s(0)=0\, ,\: \dot s(0)\equiv g$ 
are included in the definition of the balls. 
Combining the estimates of Theorems \ref{hjeqex} and \ref{u_ex}, it follows that
$$
  \| \Phi(s_1)-\Phi (s_2)\| _{\bbs ^\beta}\leq c\, T^\beta\, \| s_1-s_2\| _{\bbs ^\gb}
$$
and the contraction principle can be applied to obtain existence.
\end{proof}
\begin{rem}
It should be pointed out that this is the first well-posedness result for this class of singular 
free boundary problems in more than one space dimension for the full evolutionary problem. A 
companion result for the quasi-stationary approximation ($\varepsilon=0$) has previously 
been obtained in \cite{Gui99}.
\end{rem}

\bibliography{../../lite}
\end{document}